\documentclass{ws-jml-fixed}

\usepackage{amsmath,amssymb,amsfonts}
\usepackage{amsthm}
\usepackage{xspace,enumerate}
\usepackage{booktabs}
\usepackage{burdges-ws-jml-fixed}

\renewcommand\eq{\hbox{$^{\mathrm {eq}}$}}

\newcommand\TT{\mathcal T}
\newcommand\oo{^\circ}
\let\Union\bigcup

\renewcommand\Submitted[1][]{Submitted}

\begin{document}

\markboth{Jeff Burdges, Greg Cherlin}{Semisimple torsion in groups of finite Morley rank}

\def\version{March 2007}
\title{Semisimple torsion in groups of finite Morley rank}

\author{Jeffrey Burdges\thanks{Supported by NSF postdoctoral fellowship DMS-0503036, a Bourse Chateaubriand postdoctoral fellowship, and DFG grant Te 242/3-1.}}

\address{Institut Camille Jordan, Universit\'e de Lyon-I, Lyon, France}

\author{Gregory Cherlin\thanks{Supported by NSF grant DMS-0100794 and DMS-0600940}}

\address{Mathematics Department, Rutgers University, Piscataway, NJ, USA}

\maketitle

\begin{abstract}
We prove several results about groups of finite Morley rank without
unipotent $p$-torsion: $p$-torsion always occurs inside tori,
Sylow $p$-subgroups are conjugate, and $p$ is not the minimal
prime divisor of our approximation to the ``Weyl group.''
These results are quickly finding extensive applications within
the classification project.
\end{abstract}

\section*{Introduction}

A group of finite Morley rank is a group equipped with a notion
of dimension satisfying various natural axioms \cite[p.~57]{BN};
These groups arise naturally in model theory,
 expecially geometrical stability theory.  
The main examples are algebraic groups over algebraically closed
fields, where the notion of dimension is the usual one,
as well as certain groups arising in applications of model theory to
diophantine problems, where the notion of dimension
comes from differential algebra rather than algebraic geometry.

The main examples are algebraic groups over algebraically closed
fields, where the notion of dimension is the usual one; and the
dominant conjecture is that all such simple groups are algebraic.

\begin{namedtheorem}{Algebraicity Conjecture}[Cherlin/Zilber]
A simple group of finite Morley rank is an algebraic group
 over an algebraically closed field.
\end{namedtheorem}

Much work towards this conjecture involves local analysis
in an inductive setting reminiscent of the classification of
the finite simple groups, but without transfer arguments or
character theory.

New methods have emerged recently in the study of groups
of finite Morley rank, and have led to a number of advances.
Among the characteristic features of the recent work are the
systematic use of {\em generic covering} arguments, which
will be met with below, as well as the study of divisible abelian
$p$-subgroups (commonly known as {\em $p$-torii}), with
which we will also be occupied here.

Such $p$-tori may always be viewed as {\em semisimple}.
However, there are difficulties when one wishes to view
individual $p$-elements as either semisimple or unipotent.
For example, even a connected solvable $p$-group of a
group of finite Morley rank is merely a central product,
not necessarily a direct product, of a $p$-torus and a
definable connected nilpotent $p$-subgroups of bounded exponent,
 (commonly known as a {\em $p$-unipotent subgroup}).
Elements in the intersection have an ambiguous character.
Our main objective here is to obtain several results concerning
that $p$-torsion in connected groups of finite Morley rank
which is semisimple in a robust sense involving the absence
 of $p$-unipotent subgroups.

We say a group $G$ has {\em $p^\perp$ type} if it contains
 no nontrivial unipotent $p$-subgroup.
For the case $p=2$, one may focus on groups of $2^\perp$ type
(sometimes referred to as {\em odd or degenerate type})
because the Algebraicity Conjecture holds in the presence
of a 2-unipotent subgroup \cite{ABC}.  In this case,
our results will have numerous applications to classification
problems, beginning with the generation theorem of
\cite{BC08b}---or strictly speaking, beginning with some earlier
papers that could have been shortened had the result been
available at the time. 

\medskip

The main results are as follows.  We expect each of them to find
further use. The last three are intended to be less technical and
more readily applicable than the first two, but they do not exhaust
the information that can be extracted from the more technical results.
Two of the results stated below (Theorems 2*, 3* below) are given
in more general forms in the text.

The first section expands upon \cite{Ch05} and clarifies the
nature of the generic element of $G$. 

\begin{namedtheorem}{Theorem \ref{sec:generically_pdivisible}}
Let $G$ be a connected group of finite Morley rank, $p$ a prime, 
and let $a$ be a generic element of $G$.  Then
\begin{enumerate}
\item $a$ commutes with a unique maximal $p$-torus $T_a$ of $G$,
\item $d(a)$ contains $T_a$, and
\item If $G$ has $p^\perp$ type then $d(a)$ is $p$-divisible.
\end{enumerate}
\end{namedtheorem}

The next section contains a new genericity argument for cosets.

\begin{namedtheorem}{Theorem \ref{sec:coset_genericity}*}
Let $G$ be a group of finite Morley rank, 
let $a$ be a $p$-element in $G$ such that 
$C(a)$ has $p^\perp$ type, and
let $T$ be a maximal $p$-torus of $C(a)$ (possibly trivial).
Then 
$$ \textit{$\Union (a C\oo(a,T))^{G\oo}$ is generic in $a G\oo$.} $$ 
\end{namedtheorem}


These two technical results are the main ingredients in the
 following robust criterion for semi-simplicity.

\begin{namedtheorem}{Theorem \ref{sec:torality}*}
Let $G$ be a connected group of finite Morley rank, $p$ a prime, and
$a$ any $p$-element of $G$ such that $C_G\oo(a)$ has $p^\perp$ type.
Then $a$ belongs to a $p$-torus.
\end{namedtheorem}

Theorem \ref{sec:torality}* has immediate consequence for
 the structure of Sylow $p$-subgroups.

\begin{namedtheorem}{Corollary}
Let $G$ be a connected group of finite Morley rank of $p^\perp$ type,
and $T$ a maximal $p$-torus of $G$.
Then any $p$-element of $C(T)$ belongs to $T$.
\end{namedtheorem}

Our fourth section further exploits the genericity argument
for cosets to prove conjugacy of Sylow $p$-subgroups,
i.e.\ maximal solvable $p$-subgroups.

\begin{namedtheorem}{Theorem \ref{sec:pSylow}}
Let $G$ be a group of finite Morley rank of $p^\perp$ type.
Then all Sylow $p$-subgroups are conjugate.
\end{namedtheorem}

We note that the conjugacy of Sylow 2-subgroups in general
groups of finite Morley rank is known, and for general $p$
the result is known in solvable groups of finite Morley rank.


Our last topic concerns the so-called {\em Weyl group}, which for our
present purposes may be defined as follows.

\begin{definition}
Let $G$ be a group of finite Morley rank, and
 $T$ a maximal divisible abelian torsion subgroup of $G$.
The {\em Weyl group} of $G$ is the group $N(T)/C\oo(T)$.
\end{definition}

The maximal divisible abelian torsion subgroups of $G$ are conjugate
by \cite{Ch05}, so this group is well-defined up to conjugacy and in
particular up to isomorphism.
Furthermore it is finite since $N\oo(T)=C\oo(T)$. 

\begin{namedtheorem}{Theorem \ref{thm:Weyl_group}}
Let $G$ be a connected group of finite Morley rank.
Suppose the Weyl group is nontrivial and has odd order, with $r$ the
smallest prime divisor of its order.
Then $G$ contains a unipotent $r$-subgroup.
\end{namedtheorem}


All of these results will be needed in \cite{BC08b}, and the torality
theorem should be quite useful subsequently in the analysis of
particular configurations associated with classification problems in 
odd type groups. The corollary to Theorem \ref{sec:torality} is 
also given in \cite{BoC06} for $p=2$, and is applied there to the study of
generically multiply transitive permutation groups.

Outside material will be introduced as needed, but
 much of this occurs already in the first section.
Any facts used without explicit mention can be found in \cite{BN}.

\section{Generic $p$-divisibility}\label{sec:generically_pdivisible}

We begin by analyzing the generic element of a connected group of
finite Morley rank.  We use the notation $d(a)$ for the definable hull
of an element $a$.  The definable hull of a divisible abelian torsion
subgroup of $G$ is called a {\em decent torus}.

\begin{namedtheorem}{Theorem \ref{sec:generically_pdivisible}}
Let $G$ be a connected group of finite Morley rank, $p$ a prime, 
and let $a$ be a generic element of $G$.  Then
\begin{enumerate}
\item $a$ commutes with a unique maximal decent torus $T_a$ of $G$,
\item $d(a)$ contains $T_a$, and
\item If $G$ has $p^\perp$ type then $d(a)$ is $p$-divisible.
\end{enumerate}
\end{namedtheorem}

Here we consider only elements $a$ whose
{\em type} over $\emptyset$ is generic.
When $G$ is of $p^\perp$ type 
there need not be any generic definable {\em set} $X$ such that
 $d(a)$ is $p$-divisible for every element of $X$.
Indeed, with $G$ an algebraic torus in characteristic other than $p$,
that stronger claim fails.
In this case the generic element has infinite order,
but every infinite definable set contains 
$p$-elements of finite order, and for such elements $d(a) = \gen{a}$.

The idea of the proof is to replace the group $G$ by the centralizer
 of one of its maximal decent tori. This depends on the main result of
 \cite{Ch05}, which is closely related to point $(1)$ above.

\begin{proposition}[{\cite{Ch05}}]\label{decent_tori}
Let $G$ be a group of finite Morley rank.
Then all maximal  decent tori of $G$ are conjugate.
Furthermore, if $T$ is a maximal decent torus of $G$, then
there is a generic subset $X$ of the group $C^\o_G(T)$ such that
\begin{enumerate}
\item  $X \intersect X^g = \emptyset$ for $g\notin N_G(T)$, and
\item  $\Union X^G$ is generic in $G$.
\end{enumerate}
\end{proposition}

In particular Proposition \ref{decent_tori} states that
any element of the generic definable set $\Union_G(X^G)$ commutes with
a unique conjugate of $T$, or in other words with a unique maximal
decent torus of $G$. So we have our first point:

\begin{lemma}\label{lem:1.1}
Let $G$ be a connected group of finite Morley rank,
 and let $a \in G$ be generic over $\emptyset$.
Then $C\oo(a)$ contains a unique maximal decent torus $T_a$ of $G$.
\end{lemma}

Our next point is that $T_a$ is contained in $d(a)$, and
 at this point we must work not with the generic set $X$,
 but with the type of $a$ itself generic. 
In this situation there are two notions of genericity which
are relevant: genericity in the group $G$, and genericity in the
subgroup $C\oo(T_a)$, but again by Proposition \ref{decent_tori}
these two notions can be correlated. Indeed, the next result is a
direct reformulation of Lemma \ref{lem:1.1}.

\begin{lemma}\label{lem:1.4}
Let $G$ be a connected group of finite Morley rank,
 and $T_0$ a maximal decent torus.
Then an element $a \in C\oo(T_0)$ is generic over $\emptyset$ in $G$
 if and only if the following hold.
\begin{enumerate}
\item $T_0$ is generic over $\emptyset$, in the set of maximal decent tori;
\item $a$ is generic in the group $C\oo(T_0)$ over the canonical
  parameter for $T_0$.
\end{enumerate}
\end{lemma}


A word on terminology: as a definable set, $T_0$ can be viewed as an
``imaginary element'' of $G$, and the canonical parameter for $T_0$ is
simply this element. As this may be identified with $T_0$ itself, one
may speak of genericity ``over $T_0$''. The natural language for
discussing the group $C\oo(T_0)$ treats this parameter as a
distinguished constant; it is interdefinable with $C\oo(T_0)$.

\begin{proof}
Suppose first that $a$ is generic. 
Then $T_0=T_a$. 

We show that $T_0$ is generic over $\emptyset$.
If $T_0$ belongs to a $\emptyset$-definable family $\TT$ in $G\eq$
 (a uniformly definable family in $G$)
then $a$ belongs to the set $U := \Union_{T \in \TT}C\oo(T)$.
If the family $\TT$ is nongeneric in the set of maximal decent tori,
 then $U$ is nongeneric in $G$,
a contradiction---the failure of genericity is immediate by a rank calculation. 

Now we show $a$ is generic over the canonical parameter $t_0$
 for $T_0$, in $C(T_0)$.
If $a$ belongs to some nongeneric set $Y_{t_0} \includedin C\oo(T_0)$,
 where $Y$ is defined over $t_0$, then $a$ belongs to
 the $\emptyset$-definable nongeneric set $\Union_{t \in U} Y_{t}$,
where again the nongenericity follows by a direct rank calculation.

Now suppose $T_0$ is generic over $\emptyset$ and
 $a \in C\oo(T_0)$ is generic over the parameter $T_0$.
Suppose that $a$ belongs to the $\emptyset$-definable subset $Y$ of $G$.
Let $Y_0 = Y \intersect C\oo(T_0)$.
As $a\in Y_0$, the set $Y_0$ is generic in $C\oo(T_0)$.
The set $\TT$ of conjugates $T$ of $T_0$ for which
 $Y \intersect C\oo(T)$ is generic in $C\oo(T)$
 is $\emptyset$-definable and contains $T_0$, and
hence $\TT$ is generic in the set of conjugates of $T_0$.
It then follows from Proposition \ref{decent_tori} 
 that $Y$ is generic in $G$.
\end{proof}

We will also need some general properties of definable quotients.

\begin{lemma}\label{lem:basicdecent}
Let $G$ be a group of finite Morley rank, $A \includedin G$, 
$H$ a normal $A$-definable subgroup of $G$, and $\bar G := G/H$.
\begin{enumerate}
\item If an element $a \in G$ is generic over $A$ then its image
  $\bar a$ in $\bar G$ is generic over $A$.
\item If $T$ is a maximal decent torus of $G$, and $H$ is solvable, 
then the image $\bar T$ of $T$ in $\bar G$ is a maximal decent torus
of $\bar G$. 
\end{enumerate}
\end{lemma}

The first point has already occurred in a special form above, and is
also contained in Lemma 6.2 of \cite{Po87}.

\begin{proof}
{\em Ad 1.} Suppose $\bar V$ is an $A$-definable subset of $\bar G$
containing $\bar a$, with preimage $V$ in $G$. As $V$ contains $a$ it
is generic. But $\rk(V)=\rk(\bar V)+\rk(H)$ and thus $\bar V$ is
generic in $\bar G$. Thus $\bar a$ is generic in $\bar G$.

{\em Ad 2.} 
Let $T_p$ be the $p$-torsion subgroup of $T$.
It suffices to show that $\bar T_p$ is a maximal $p$-torus of $\bar G$.
Let $S_p$ be the preimage in $G$ of a maximal $p$-torus $\bar S_p$
 of $\bar G$ containing $\bar T_p$.
We may suppose that $G=d(S_p)$ and thus $G$ is solvable. 
Now let $P$ be a Sylow $p$-subgroup of $G$ containing $T_p$.
Then $\bar P$ is a Sylow $p$-subgroup of $\bar G$ \cite{ACCN98}
 and therefore contains a maximal $p$-torus of $\bar G$. 
Hence $\bar P$ contains a maximal $p$-torus of $\bar G$.
But as $S$ is a solvable $p$-group,
 $S\oo = T_p*U_p$ with $U_p$ unipotent \cite[Corollary 6.20]{BN},
 so $\bar T_p$ is the maximal $p$-torus of $\bar P$,
and hence is a maximal $p$-torus of $\bar G$.
\end{proof}


\begin{lemma}\label{lem:1.2}
Let $G$ be a connected group of finite Morley rank
and $a \in G$ generic.  Then $d(a)$ contains $T_a$.
\end{lemma}

\begin{proof}
Treating the parameter $T_a$ as a constant, 
and bearing in mind Lemma \ref{lem:1.4},
we may suppose that $a$ is
a generic element of $C\oo(T_a)$, and $T_a$ is $\emptyset$-definable.
Hence we may replace $G$ by
$C\oo(T_a)$, assuming therefore that 
$$\hbox{$G$ contains a unique maximal decent torus $T$,
  which is central in $G$.}$$ 

Let $T_1$ be the definable hull of
 the torsion subgroup of $d(a) \intersect T$.
As $T$ is taken to be $\emptyset$-definable,
 the torsion subgroup of $T$ is contained in $\acl(\emptyset)$
and hence the definable set $T_1$, treated as another parameter,
 is also in $\acl(\emptyset)$.
Thus $\bar a$ is generic in the quotient $\bar G := G/T_1$,
 and in this quotient $\bar T := T/T_1$ is a maximal decent torus.
So replacing $G$ by $\bar G$, we may suppose that $d(a)\intersect T$
 is torsion free. It suffices to show that $T=1$.

By \cite[Ex.~10 p.~93]{BN}, $d(a)=A \oplus C$ is the direct
sum of a divisible abelian group $A$ and a finite cyclic group $C$. 
If $n=|C|$ then for any multiple $N$ of $n$, we have $d(a^N)=A$.
On the other hand, for any torsion element $c \in T$,
 the element $a'=ac$ is also generic over $\emptyset$
and hence $a'$ and $a$ realize the same type.
Letting $N$ be a multiple of $n$ and the order of $c$,
 it follows that $d((a')^n)=d((a')^N)=d(a^N)=d(a^n)$
 and thus $c^n \in d(a^n) \le d(a)$.
Now by our reductions $d(a)$ contains no $p$-torus for any $p$,
 and hence $d(a)$ has bounded exponent. 
Thus $c^n$ has bounded exponent, with $c$ varying and $n$ fixed,
 and so $T=1$ as claimed.
\end{proof}

For the final point in Theorem \ref{sec:generically_pdivisible}
we prepare the following, which is a minor variation on a result 
of \cite{BBC}.

\begin{lemma}\label{central_maxtorus}
Let $G$ be a connected group of finite Morley rank, $p$ a prime, 
and $T$ a maximal $p$-torus of $G$.
Suppose that $T$ is central in $G$ and
 $a$ is a $p$-element of $G$ not in $T$.
Then $C\oo(a)$ contains a nontrivial $p$-unipotent subgroup.
Thus if $G$ is of $p^\perp$ type then all $p$-elements in $C(T)$
belong to $T$.
\end{lemma}

Here we employ one of the main results of \cite{BBC}.

\begin{fact}[{\cite[Theorem 4]{BBC}}]\label{BBC_p}
Let $G$ be a connected group of finite Morley rank, and
 let $a \in G$ be a $p$-element.
Then $C(a)$ contains an infinite abelian $p$-subgroup.
\end{fact}

\begin{proof}[Proof of Lemma \ref{central_maxtorus}]
Observe first that the $p$-torsion subgroup of $d(T)$ is $T$, and thus
$a \notin d(T)$. Now passing to a quotient as in the previous argument
we may suppose that $T=1$ and $G$ contains no $p$-torus.
So $C\oo(a)$ contains a nontrivial $p$-unipotent subgroup
 by Fact \ref{BBC_p}.
\end{proof}

We turn to the last point in Theorem
\ref{sec:generically_pdivisible}.

\begin{lemma}\label{lem:1.3}
Let $G$ be a connected group of finite Morley rank of $p^\perp$ type, 
and $a \in G$ a generic element.
Then $d(a)$ is $p$-divisible.
\end{lemma}

\begin{proof}
As we have seen above, we may suppose that $T_a$ is central in $G$
 and $\emptyset$-definable.
The group $d(a)$ is an abelian group of finite Morley rank,
and hence has the form $A \times C$ for
 some $p$-divisible abelian group $A$, and
 some $p$-group $C$ of bounded exponent by \cite[Ex.~10, p.~93]{BN}. 
Since $G$ is of  $p^\perp$ type,
 $C \leq T_a$ by Lemma \ref{central_maxtorus}.
As $T_a \leq d(a)$ by Lemma \ref{lem:1.2},
 $T_a \leq A$ and $d(a) = A$ is $p$-divisible.
\end{proof}

Now Theorem \ref{sec:generically_pdivisible} is contained in 
Lemmas \ref{lem:1.1}, \ref{lem:1.2}, and \ref{lem:1.3}.

\section{Coset genericity}\label{sec:coset_genericity}

In this section we prove a {\em generic covering theorem}.
Theorems of this type have played an increasing role in the
analysis of connected groups of finite Morley rank.
Our aim here is to show that for a $p$-element $a$ of a group
$G$ of $p^\perp$ type, the union of the conjugates of $C\oo(a)$
is generic in $G$.  This improves on the analysis carried out in
\cite{BBC} for groups of $p$-degenerate type. 
In order to prove this we need to sharpen it substantially and
identify a subgroup of $C\oo(a)$ actually responsible for the
genericity. The precise result we aim at is the following, which
generalizes the result in several directions, notably by allowing
the element $a$ to lie outside the connected component of $G$.

We formulate this analysis using a more general set of primes $\pi$,
as opposed to the single prime $p$ used in the introductory
statement.  A $\pi$-torus is a divisible abelian $\pi$-group.
Similarly $\pi^\perp$ type means $p^\perp$ type for all $p \in \pi$.

\begin{namedtheorem}{Theorem \ref{sec:coset_genericity}}
Let $G$ be a group of finite Morley rank, 
let $a$ be a $\pi$-element in $G$ such that 
$C(a)$ has $\pi^\perp$ type, and
let $T$ be a maximal $\pi$-torus of $C(a)$ (possibly trivial).
Then 
$$ \textit{$\Union (a C\oo(a,T))^{G\oo}$ is generic in $a G\oo$.} $$ 
In particular
$$\textit{$\Union (a C\oo(a)^{G\oo})$ is generic in $a G\oo$.}$$ 
\end{namedtheorem}

Generic covering theorems have involved definable subgroups more
often than cosets. The following covering lemma, given in \cite{BBC},
is well adapted to the case of cosets.

\begin{fact}[{\cite[Lemma 4.1]{BBC}}]\label{generic_covering}
Let $G$ be a group of finite Morley rank, $H$ a definable subgroup of $G$,
 and $X$ a definable subset of $G$. Suppose that
$$\rk(X\setminus \bigcup_{g \notin H} X^g) \geq \rk(H)$$
Then $\rk(\Union X^G) = \rk(G)$.
\end{fact}

The following property of generic subsets of cosets is very well
known for subgroups, but occurs more rarely in its general form.

\begin{lemma}
Let $G$ be a group of finite Morley rank, $H$ a connected definable
subgroup, and $X$ a definable generic subset of the coset $a H$. 
Then $\gen{X}=\gen{a H}=\gen{a,H}$.
\end{lemma}

\begin{proof}
The second equality is purely algebraic, and clear.
For the first, an application of genericity and connectedness shows
that $H \le \gen{X}$, and thus $a H\includedin \gen{X}$.
\end{proof}

\begin{proof}[Proof of Theorem \ref{sec:coset_genericity}]
We will use the notation $N_G(X)$ here for arbitrary subsets of $G$,
not just subgroups, with its usual meaning: the setwise stabilizer of
$X$ under the action of $G$ by conjugation.

Let $\TT$ be the set of maximal $\pi$-tori of $C\oo_G(a)$.
We observe first that $\TT$ may be identified with 
a definable set in $G\eq$. Indeed, it follows from the conjugacy of
maximal decent tori that maximal $\pi$-tori are conjugate under the
action of the group 
$$G_a=C\oo_G(a)$$
so $\TT$ corresponds naturally to the right coset space
$N_{G_a}(T)\backslash G_a$ for any $T\in \TT$, and $N(T)=N(d(T))$ is definable.
As the elements of $\TT$ are themselves undefinable, this
identification should be used with circumspection.  

As the maximal $\pi$-tori of $C(a)$ are conjugate, we may suppose that
the $\pi$-torus $T \in \TT$ is chosen generic over $a$. 
Set $$ H := C\oo(\gen{a,T}) $$
which enters the picture most naturally here as 
$C_{C\oo(a)}\oo(T)$.
Then $T$ is the unique maximal $\pi$-torus in $H$,
and we aim to show that 
$$\rk(\Union (a H)^{G\oo})=\rk(G\oo)$$

Let $\hat H$ be the generic stabilizer of $a H$, 
defined as
$$\{g \in G: \rk((a H) \intersect (a H)^g) = \rk(a H)\}$$
This is a definable subgroup of $G$.
We claim 
$$\rk(\hat H) = \rk(H)\leqno(1)$$

Since $a$ is an element of finite order
normalizing (even centralizing) $H$, the group 
$\gen{a,H}$ is definable with $\gen{a,H}\oo =H$.
Applying the preceding lemma,
$$\hat H \leq N_G(\gen{a,H}) \le N_G(\gen{a,H}\oo)=N_G(H) \leq N(T)$$
Thus $\hat H\oo \le C(T)$.

We claim that any $\pi$-element $u$ of $\gen{a,H}$ 
lies in the abelian group $\gen{a,T}$:
indeed, the $\pi$-group $U = \gen{u,a}$ has the form
$U_0 \gen{a}$ with $U_0=U \intersect H$.
We claim that $U_0 \le T$.
For this, it suffices to show that any $\pi$-element $u' \in U_0$
belongs to $T$.  But this holds by Lemma \ref{central_maxtorus}. 

Therefore $\gen{a,H}$ contains only finitely many elements of the same
 order as $a$, and as $\hat H$ acts by conjugation on these elements, 
we have $\hat H\oo \leq C\oo(a)$
 and thus $\hat H\oo \leq C\oo(\gen{a,T}) = H$.
So $(1)$ holds.

We would like to 
apply the generic covering lemma, Fact \ref{generic_covering}, 
with $X = a H$ and with $H$ (in the lemma) equal to $\hat H$ (here).
For this, it suffices to verify the condition
$$\rk(a H \setminus \Union_{g \notin \hat H} (a H)^g) = \rk(H) \leqno(*)$$

Now suppose $x \in a H$ is generic over the parameters $a$ and $T$
(really, $d(T)$).
We claim that both
$$\hbox{$x$ centralizes a unique maximal $\pi$-torus of $C\oo(a)$, namely $T$, and}\leqno(2)$$
$$a\in d(x)\leqno(3)$$

Clearly $x = a h$ with $h \in H$ generic over $a$ and $T$.
Since $T$ is itself generic over $a$, $h$ realizes the type of
 a generic element of $C\oo(a)$ over $a$ (Lemma \ref{lem:1.4}).
By Theorem \ref{sec:generically_pdivisible} for $p \in \pi$,
 $h$ centralizes a unique maximal $\pi$-torus of $C\oo(a)$ for $p \in \pi$,
 and hence so does $x$. So $(2)$ follows.

As $C\oo_G(a)$ has $\pi^\perp$ type,
 $d(h)$ is $\pi$-divisible by Theorem \ref{sec:generically_pdivisible}.
So, for any $\pi$-number $q$,
 the quotient $d(h)/d(h^q)$ is a $p$-divisible group of exponent at most $q$,
 and hence trivial: $d(h)=d(h^q)$.
Let $q$ be the order of $a$.
Then $d(x^q) = d(a^q h^q) = d(h^q) = d(h)$.
So $h \in d(x)$, and hence also $a \in d(x)$, giving $(3)$.

If $(*)$ fails, then
 $H\intersect \Union_{g \notin \hat H}(a H)^g$ is generic in $a H$, 
so, as $x \in a H$ is generic over the parameters $a$ and $T$,
 we have  $x^g \in a H$ for some $g \notin \hat H$.
As $x,x^g \in a H$, $d(x)$ and $d(x^g)$ both commute with $T$,
 and therefore $d(x)$ also commutes with $T^{g^{-1}}$.
Since $a \in d(x)$, we have $T^{g^{-1}} \leq C\oo_G(a)$.
By $(2)$ it follows that $T=T^{g^{-1}}$, that is $g \in N(T)$.

Again, since $a \in d(x)$,
 we have $a^g \in d(x^g)$ is an element of order $q$ in $\gen{a,H}$.
So $a^g \in \gen{a,T}$.
Since $g \in N(T)$ this gives $g \in N_G(\gen{a,T})$ as well.
Now $H=C\oo(\gen{a,T})$ so 
$$H^g=C\oo(\gen{a,T})^g=C\oo(\gen{a^g,T})=C\oo(\gen{a,T})=H$$

Hence $(a H)^g = (xH)^g = x^g H^g = x^g H = a H$, and $g \in \hat H$,
a contradiction.
So $(*)$ holds, and our result follows by Fact \ref{generic_covering}. 
\end{proof}

\section{Torality}\label{sec:torality}

We now prove the main result of the paper.
Again, we formulate this in a technical form slightly more general
than the original statement, using a set of primes $\pi$.

\begin{namedtheorem}{Theorem \ref{sec:torality}}
Let $G$ be a connected group of finite Morley rank, $\pi$ a set of primes,
and $a$ any $\pi$-element of $G$ such that $C_G\oo(a)$ has $\pi^\perp$ type.
Then $a$ belongs to a $\pi$-torus.
\end{namedtheorem}

Theorem \ref{sec:torality} has the following direct corollary.

\begin{corollary}\label{cor:torality}
Let $G$ be a connected group of finite Morley rank with
 a $\pi$-element $a$ such that $C(a)$ has $\pi^\perp$ type.
Then $a$ belongs to any maximal $\pi$-torus of $C(a)$.
\end{corollary}

\begin{proof}
By Theorem \ref{sec:torality}, there is a $\pi$-torus $T$ containing $a$.
By Proposition \ref{decent_tori}, 
 any maximal $\pi$-torus in $C(a)$ is $C\oo(a)$-conjugate to $T$,
and so contains $a$.
\end{proof}

This imposes very strong restrictions on
 a simple group $G$ of finite Morley rank.
For $\pi=\{2\}$, the outstanding structural problems
 concern groups of $2^\perp$ type (i.e., odd or degenerate type).
In this context, our results impose constraints on the
 structure of a Sylow 2-subgroup,
 which will be developed in \cite{BC08b}.

For the proof, we use the following variation
 on Fact \ref{BBC_p} \cite[Theorem 4]{BBC}.
This lemma is due to Tuna \Altinel.

\begin{lemma}[\Altinel]\label{BBC_pi}
Let $G$ be a connected group of finite Morley rank, and
 let $a \in G$ be a $\pi$-element.
Then $C(a)$ contains an infinite abelian $p$-subgroup for some $p \in \pi$.
\end{lemma}

We require the following.

\begin{fact}[{\cite{ABCCdraft}; \cite[\qFact 3.2]{Bu03}}]\label{Cquotient}
Let $G = H \rtimes T$ be a group of finite Morley rank
 with $H$ and $T$ definable.
Suppose $T$ is a solvable $\pi$-group of bounded exponent and
$Q \normal H$ is a definable solvable $T$-invariant $\pi^\perp$-subgroup.
Then $$ C_H(T)Q/Q = C_{H/Q}(T)\mathperiod $$
\end{fact}


\begin{proof}[Proof of Lemma \ref{BBC_pi}]
We may take $G$ to be a minimal counterexample;
 in particular $C\oo(a)$ is a $\pi^\perp$-group.
Of course, $G$ does contains an infinite abelian $p$-group
 for some $p \in \pi$ by Fact \ref{BBC_p}.
So clearly $a \notin Z(G)$.

As $Z\oo(G)$ has no $\pi$-torsion,
 $C_{G/Z\oo(G)}(a) = C_G(a) / Z\oo(G)$ by Fact \ref{Cquotient}.
So $C_{G/Z\oo(G)}(a)$ has no $\pi$-torsion
 by \cite[Ex.~11 p.~93 or Ex.~13c p.~72]{BN}. 
Thus $Z\oo(G) = 1$ by minimality of $G$.

We now show that $a \in d(x) \cap a C\oo(a)$ for any $x \in a C\oo(a)$.
Let $K := d(x) \cap C\oo(a)$.  So $x$ is a $\pi$-element in $d(x)/K$.
By \cite[Ex.~11 p.~93]{BN}, 
 $x d\oo(x)$ contains a $\pi$-element $b$.
But $a$ is the unique $\pi$-element in $a C\oo(a) \supseteq x K$.
Thus $a = b \in d(x)$, as desired.

By Theorem \ref{sec:coset_genericity},
 $\bigcup (a C\oo(a))^G$ is generic in $G$.
We show that $G$ has no divisible torsion.
Otherwise, choose a maximal decent torus $T$ of $G$.
By Fact \ref{decent_tori},
 $\bigcup C\oo(T)^G$ is generic in $G$ too, and
 hence meets $a C\oo(a)$ in an element $x$.
So $a \in d(x)$ lies inside some $C\oo(T)^g$.
But $C_{C\oo(T)^g}(a)$ is still a $\pi^\perp$-group,
 contradicting the minimality of $G$.

As $C\oo(a^{-1}) = C\oo(a)$,
 $\bigcup (a^{-1} C\oo(a))^G$ is also generic in $G$,
 by Theorem \ref{sec:coset_genericity}.
So there is some $x \in a^{-1} C\oo(a) \cap (a C\oo(a))^g$ for some $g \in G$.
As above $a^g$ and $a^{-1}$ are the only $\pi$-elements
 in $(a C\oo(a))^g$ and $a^{-1} C\oo(a)$, respectively.
So $a^g = a^{-1}$.
It follows that $G$ has an involution in $d(g)$.

We recall that $B(G)$ is the subgroup of $G$ generated by
 all its 2-unipotent subgroups.
As $G$ has no divisible torsion,
 $G$ has even type, but has no algebraic simple section.
So $B(G)$ is a 2-unipotent subgroup normal in $G$,
 by the Even Type Theorem \cite{ABC}.
Now $Z\oo(B(G)) \neq 1$ by cite[Lemma 6.2]{BN}.

As $G$ has no divisible torsion, 
 it is generated by $p$-unipotent and $(0,r)$-unipotent
 nilpotent subgroups, for various $p$ and $r$.
All these nilpotent $(0,r)$-unipotent subgroups
 centralize $B(G)$ by \cite[Lemma 4.5]{Bu05a}.
For $p \neq 2$, all nilpotent $p$-unipotent subgroups
 centralize $B(G)$ too.
Therefore $Z\oo(B(G)) \leq Z\oo(G) = 1$, a contradiction.
\end{proof}


\begin{proof}[Proof of Theorem \ref{sec:torality}]
By Lemma \ref{BBC_pi},
 there is a non-trivial $\pi$-torus $T$ of $C\oo(a)$,
 which we take maximal in $C\oo(a)$.  Set $H := C\oo(a,T)$.
By Theorem \ref{sec:coset_genericity},
 the set $(a H)^G$ is generic in $G$.
So after conjugating
 we may suppose that some $x = a h \in aH$ is generic in $G$.
We claim that $a \in d(x)$.

Since $x$ is generic in $G$, 
 $C(x)$ contains a unique maximal $\pi$-torus $S$ of $G$,
 which lies inside $d(x)$,
 by Theorem \ref{sec:generically_pdivisible}.
Clearly $T \leq S$ since $T \le C(x)$.
The definable hull $d(x)$ contains a $\pi$-element $x'$
 with $x' H = x H = a H$.
So $x' a^{-1} \in H$ is a $\pi$-element, and lies inside $C(T)$.
Thus $x' a^{-1} \in S \le d(x)$ by Lemma \ref{central_maxtorus},
 and so $a \in d(x)$.

Again since $x$ is generic in $G$,
 we have $x \in C\oo(S)$ by Theorem \ref{sec:generically_pdivisible},
and hence $a \in C\oo(S)$ and $T = S$ is nontrivial.
Since $H \leq C(a)$ has $\pi^\perp$ type,
 we have $a \in T$ by Lemma \ref{central_maxtorus}.
This proves our claim.
\end{proof}

For applications to the structure of Sylow $p$-subgroups in
 connected groups of $p^\perp$ type and low Pr\"ufer $p$-rank,
especially Pr\"ufer rank 1, see \cite{BC08b}.

\section{Conjugacy of Sylow $p$-subgroups}\label{sec:pSylow}

We define {\em Sylow $p$-subgroups} of a group $G$ of finite
Morley rank as maximal solvable $p$-subgroups.
One arrives at the same class of groups by imposing local finiteness
or local nilpotence in place of solvability \cite[\S6.4]{BN}.
If $S$ is a  Sylow $p$-subgroup of $G$ then $S\oo$ will be a
central product of a $p$-unipotent subgroup and a $p$-torus,
and in particular is nilpotent.
So if $S$ is a Sylow $p$-subgroup and $X$ a proper subgroup of
$S$, then $N_S(X)>X$. 

Our goal in the present section is the following.

\begin{namedtheorem}{Theorem \ref{sec:pSylow}}
Let $G$ be a group of finite Morley rank of $p^\perp$ type.
Then all Sylow $p$-subgroups are conjugate.
\end{namedtheorem}

The conjugacy result is also known for solvable groups, as a special
case of the theory of Hall subgroups (\cite[Theorem 9.35]{BN})
and for arbitrary groups of finite Morley rank when $p=2$. 

As an immediate consequence we can strengthen
\cite[Theorem 3]{BBC}.

\begin{corollary}\label{BBC_p_plus}
Let $G$ be a connected group of finite Morley rank and $p^\perp$ type.
If $G$ some Sylow $p$-subgroup of $G$ is finite
then $G$ contains no elements of order $p$. 
\end{corollary}

The critical case for the proof is the case in which at least one
Sylow $p$-subgroup is finite; which proves the corollary itself.
It also shows that Sylow $p$-subgroups are conjugate if all lie outside $G\oo$.

\begin{lemma}
Suppose $G$ is a group of finite Morley rank and $p^\perp$ type
containing a finite Sylow $p$-subgroup $P$.
Then all Sylow $p$-subgroups of $G$ are conjugate.
\end{lemma}

\begin{proof}
Let $O_p(G)$ denote the subgroup of $G$ generated by its solvable
normal $p$-subgroups. Such a $p$-subgroup must be contained in $P$ and
thus $O_p(G)\le P$ is finite, and is the largest finite normal
$p$-subgroup of $G$. In $\bar G=G/O_p(G)$ we have $O_p(\bar G)=1$ and
$\bar P=P/O_p(G)$ is a finite Sylow $p$-subgroup of $\bar G$, and if
we prove the claim for $\bar G$ it follows for $G$. So we may suppose
$$O_p(G)=1\leqno(1)$$

Let $D$ be a subgroup of $P$ of maximal order subject to the 
condition: $D$ is contained in a solvable $p$-subgroup of $G$ which has no
conjugate contained in $P$. Let $R$ be such a $p$-subgroup. Let
$D_1=N_P(D)$, $D_2=N_R(D)$. By the maximality of $D$, any $p$-Sylow
subgroup $P_1$ of $N(D)$ containing $D_1$ is conjugate to a subgroup of
$P$. Let $R_1$ be a Sylow $p$-subgroup of $N(D)$ containing $D_2$. If
$R_1$ is conjugate to $P_1$, then $R_1$ is conjugate to a subgroup of
$P$. In particular $D_2$ is then conjugate to a subgroup of $P$ and
$R$ is conjugate to a group meetings $P$ in a subgroup of order
greater than $|D|$. But this contradicts the choice of $D$. 

It follows that in $N(D)$ we have nonconjugate Sylow $p$-subgroups, so
by the minimality of $G$ we find $D\normal G$ and thus $D\le O_p(G)=1$.
Hence any solvable $p$-subgroup which meets $P$ nontrivially is
conjugate to a subgroup of $P$.

Fix $a \in P$ nontrivial. We claim
$$\hbox{$C\oo(a)$ is a $p^\perp$-group}\leqno(2)$$
If this fails, take $x \in C\oo(a)$ a nontrivial $p$-element. 
By Fact \ref{BBC_p},  
 $C(x)$ contains an infinite abelian $p$-subgroup $A$.
As $O_p(G)=1$, we have $C(x)<G$ and hence the
Sylow $p$-subgroups of $C(x)$ are conjugate.
Taking Sylow $p$-subgroups  $Q$ and $R$ of $C(x)$
 containing $\gen{a,x}$ and $A$ respectively,
we find that $Q$ is conjugate to a subgroup of $P$
 since $Q$ meets $P$ nontrivially, and
hence the infinite group $R$ is conjugate to a
 subgroup of the finite group $P$, a contradiction. 

Now let $b$ be an arbitrary $p$-element of the coset $a G\oo$,
 and $T_b$ a maximal torus of $C\oo(b)$.
Then $\Union (b C\oo(b,T_b))^{G\oo}$ is generic in $a G\oo$
 by Theorem \ref{sec:coset_genericity}.
This applies in particular to $a$, with $T_a=1$.
As we have generic sets associated
 to $a$ and $b$ in the coset $a G\oo$,
their intersection is nontrivial, giving
$$ a C\oo(a)\intersect b'C\oo(b',T_{b'})>1
\quad\textrm{for some conjugate $b'$ of $b$.} $$

Fix $h \in aC\oo(a)\intersect b'C\oo(b',T_{b'})$. 
Since $h \in a G\oo$, there is a $p$-element $h'\in d(h)\intersect
a G\oo$. But 
$d(h)\intersect a G\oo$ is contained in both
$aC\oo(a)$ and $b'C\oo(b',T_{b'})$.
So $h'\in aC\oo(a)$, and as $h'$ is a $p$-element
we find $h'=a$. Similarly $h'\in b'T_{b'}$. Thus $a\in b'T_{b'}$ 
and $T_{b'}\le C\oo(a)$, so $T_{b'}=1$. We conclude that $a=b'$ and
thus 
$$\hbox{For $a \in P^\#$, any two $p$-elements in $a G\oo$ are conjugate.}
\leqno(3)
$$

Now fix an arbitrary Sylow $p$-subgroup $Q$ of $G$. 
We will show that $P$ and $Q$ are conjugate.

Let $\bar G=G/G\oo$ and let $\bar R$ be a Sylow $p$-subgroup of $\bar
R$ containing $\bar P$. We may suppose after conjugating $Q$ that
$\bar Q\le \bar R$. We claim
$$\bar R=\bar P\leqno(4)$$
Assuming the contrary, let $R$ be the preimage in $G$ of $\bar R$. 
We have $N_{\bar R}(\bar P)>\bar P$ and thus $N_R(PG\oo)>PG\oo$. 
If $PG\oo<G$ then Sylow $p$-subgroups of $PG\oo$ are conjugate and 
therefore $N_R(PG\oo)=G\oo N_R(P)$. Thus $N_R(P)$ covers $N_{\bar
  R}(\bar P)$ and therefore there is a $p$-element $x\in
N_R(P)\setminus P$. But then $P$ is not a Sylow $p$-subgroup, a
contradiction. So (4) holds.

Hence $PG\oo=QG\oo$. In particular there are $a\in P^\#$, $b\in Q$ with 
$a G\oo=bG\oo$ and thus $a,b$ are conjugate. Hence some conjugate of
$Q$ meets $P$, and as we have shown this conjugate of $Q$ must itself
be conjugate to $P$.
\end{proof}

\begin{proof}[Proof of Theorem \ref{sec:pSylow}]
We have $G$ a group of finite Morley rank of $p^\perp$ type and $P_1,P_2$ 
Sylow $p$-subgroups. We may suppose that Sylow $p$-subgroups in proper
definable subgroups of $G$ are conjugate, and we wish to prove the
same for $G$.

Let $T_1,T_2$ be the maximal $p$-tori in
$P_1,P_2$ respectively. We may suppose $T_1\le T_2$.
If $P_1$ is finite the preceding lemma applies. So we may suppose
$T_1$ is nontrivial. 

If $N_G(T_1)<G$ then conjugacy holds in $N_G(T_1)$ and thus $T_2$ is
conjugate to a subgroup of $P_1$. In this case $T_1=T_2$, so
$P_1,P_2\le N(T_1)$ and our claim follows.

So suppose $T_1\normal G$. Then passing to $\bar G=G/d(T_1)$, the
image $\bar P_1$ of $P_1$ is finite. We claim that $\bar P_1$ is a
Sylow $p$-subgroup of $\bar G$. Let $\bar Q_1$ be a solvable $p$-group
containing $\bar P_1$, set $\bar Q=d(\bar Q_1)$, and let $Q$ be the
preimage in $G$ of $\bar Q$. Then $Q$ is solvable. By \cite{ACCN98}, 
$\bar P_1$ is a Sylow $p$-subgroup of $\bar Q$, and hence $\bar
P_1=\bar Q_1$. That is, $\bar P_1$ is a Sylow $p$-subgroup of $\bar
G$. 

By the previous lemma $\bar P_1$ and $\bar P_2$ are conjugate, and we
may suppose they are equal. Let $\bar P=d(\bar P_1)$ and let $P$ be
the preimage in $G$ of $\bar P$. Then $P$ is solvable and $P_1,P_2\le
P$, so by \cite{ACCN98} the groups $P_1,P_2$ are conjugate, as claimed.
\end{proof}

\section{Weyl groups}\label{sec:Weyl_group}

A suitable notion of ``Weyl group'' in the context of groups of finite
Morley rank is the following.

\begin{definition}
Let $G$ be a group of finite Morley rank.
Then the {\em Weyl group} associated to $G$ is
 the abstract group $W = N(T)/C\oo(T)$
where $T$ is a maximal decent torus.
\end{definition}

This is well-defined up to conjugacy in $G$, and finite.

In algebraic groups, Weyl groups are Coxeter groups, and in particular
are generated by involutions.
We note that, by a Frattini argument (Proposition \ref{decent_tori}),
 the ``Weyl group'' associated to some non-maximal decent torus
 is a quotient of the Weyl group associated to a maximal decent torus,

\begin{theorem}\label{thm:Weyl_group}
Let $G$ be a connected group of finite Morley rank.
Suppose the Weyl group is nontrivial and has odd order,
 with $r$ the smallest prime divisor of its order.
Then $G$ contains a unipotent $r$-subgroup.
\end{theorem}

In fact, we prove that either
\begin{enumerate}
\item[(H1)] Any $r$-element representing an element of order $r$ in $W$
 centralizes a unipotent $r$-subgroup,
or else
\item[(H2)] Some toral $r$-element centralizes a unipotent $r$-subgroup.
\end{enumerate}


\begin{corollary}
Let $G$ be a minimal connected simple group of finite Morley rank.
Suppose the Weyl group is nontrivial and has odd order,
 with $r$ the smallest prime divisor of its order.
Then (H1) holds in $G$.
\end{corollary}

\begin{proof}
Otherwise there is an $r$-element $x$ representing
 a Weyl group element of order $r$ such that $C(x)$ has $p^\perp$ type.
Then $x$ is toral by Theorem \ref{sec:torality}.
By (H2), there is a toral $r$-element $y$ such that $U_p(C(y)) \neq 1$.
We may assume that $[x,y] = 1$ by \cite{Ch05}.
But now $x \in C\oo(y)$.  So $Z(U_p(C(y))) \leq C(x)$, a contradiction.
\end{proof}

Here one need only assume that simple sections have $p^\perp$ type,
 not full minimal simmplicity.

\begin{proof}[Proof of Theorem \ref{thm:Weyl_group}]
We consider a counterexample $G$ with maximal decent torus $T$, 
set $W := N(T)/C\oo(T)$, and take $r | \abs{W}$ minimal.
Suppose also that both (H1) and (H2) fail in $G$.
We take $G$ to have minimal Morley rank subject to these conditions.

We first reduce to the case
$$Z(G)=1$$
In $\bar G = G / Z(G)$ the image $\bar T$ of $T$ is
 a maximal decent torus by Lemma \ref{lem:basicdecent},
 with preimage $T Z(G)$ in $G$, and
$T$ is the unique maximal decent torus of $T Z(G)$.
Hence $N(\bar T)$ is the image of $N(T)$ and so
 $C\oo(\bar T) = N\oo(\bar T)$ is the image of $C\oo(T) = N\oo(T)$.
Thus $N(\bar T) / C\oo(\bar T)$ is isomorphic with $W$,
So we may assume $Z(G) = 1$ after replacing $G$ by $\bar G$.

Now let $T_r$ be the maximal $r$-torus of $T$,
 which is nontrivial by \cite[Theorem 3]{BBC}.
Fix an element $w$ of order $r$ in $W$, and
choose a representative $a \in N(T)$ for $w$
 which is itself an $r$-element.
We now assume that $C\oo(a)$ has $r^\perp$ type since (H1) fails.
Then $a \in N(T_r) \setminus T_r$.

We claim
$$\hbox{$C_{T_r}(a)$ is finite.}$$

Otherwise, set $U := C\oo_{T_r}(a) \neq 1$.
Let $\hat U$ be a maximal $r$-torus of $C\oo_G(a)$
 which contains $U$.
Then $a \in \hat U$ by Theorem \ref{sec:torality}.
So $H = C\oo_G(U)$ contains both $T$ and $a$, since $a \in \hat U$.
As $Z(G)=1$, we have $H < G$.
The Weyl group of $H$ is $N_H(T)/C_H(T)$, and
 $a$ represents an $r$-element of this group,
This contradicts the supposed minimality of $G$.
So indeed $C_{T_r}(a)$ is finite.

Now commutation with $a$ produces an endomorphism of $T_r$
 with finite kernel, and
it is easy to see that any such endomorphism is surjective
(working either which large invariant finite subgroups of $T_r$, or
equivalently with the action of the endomorphism ring of $T_r$ on
the dual ``Tate module'').
So $[a,T_r]=T_r$, and multiplying by $a$ on the left gives
$$ a^{T_r}=a\cdot T_r \leqno(*) $$
Now there is some $r$-element $b \in C_{T_r}(a)^\#$ since
$\Omega_1(T_p)\cdot \gen{a}$ is a finite $p$-group.
Our goal is to play with $(*)$ and variations of $(*)$ to show that
$b$ and $b^2$ are conjugate,
 under the action of the Weyl group,
which will contradict our hypothesis on the minimality of $r$. 
Observe to begin with that $a$, $a b$, and $a b^2$ 
are all $T_r$-conjugate, as they are in the coset $a T_r$.

We show next that $a \notin C\oo(b)$ and $b \notin C\oo(a)$.

If $a \in C\oo(b)$ then, as $T \le C\oo(b)$,
 we again contradict the minimality of $G$.
Hence $a \notin C\oo(b)$.
On the other hand we may now assume, by failure of (H2),
 that $C\oo(b)$ has $r^\perp$ type since $b$ is toral.
So, if $b \in C\oo(a)$ then,
 by Theorem \ref{sec:torality} and its corollary,
$b$ belongs to a maximal $r$-torus $U$ of
$C\oo(a)$, and also $a \in U$; so $a \in C\oo(b)$, a contradiction.
Thus $b \notin C\oo(a)$.

As $C(a)$ has $r^\perp$ type, Theorem \ref{sec:pSylow}
 says that its Sylow $r$-subgroups are conjugate.
So, as $b$ lies inside $C(a)$,
 there is a maximal decent torus $U$ of $C(a)$ normalized by $b$.
Now $a \in U$ by Corollary \ref{cor:torality},
 and $U$ is a maximal decent torus of $G$.
Thus $b$ represents a nontrivial $r$-element of
 the Weyl group relative to $U$.
Now, with $U_r$ the $r$-torsion in $U$,
 we can reverse the roles of $a$ and $b$,
and conclude that $b,ab$ are conjugate under the action of $U_r$.
Thus $a,b$ are conjugate in $G$.
As $r>2$ we may argue similarly that $a,b^2$ are conjugate in $G$.
So $b,b^2$ are conjugate in $G$.

Now observe, by \cite[Lemma 10.22]{BN}, that
 $N(T)$ controls fusion in $T$:
if $X \includedin T$ and $X^g \includedin T$ then
 $T,T^g \le C(X^g)$ and hence $T^g$ is conjugate to $T$ in $C(X^g)$,
thus $T$ is conjugate to $T^g$ in $N(T)$.
So $b$ and $b^2$ are conjugate under the action of $N(T)$.
In other words, we have a Weyl group element carrying $b$ to $b^2$.
Thus we have elements in $W$ whose order is some prime dividing
 the order of 2 in the multiplicative group modulo $r$.
Such a prime is a factor of $r-1$, and hence less than $r$.
This contradicts the minimization of $r$.
\end{proof}

\begin{corollary}\label{cor:Weyl_group}
Let $G$ be a connected group of finite Morley rank 
 without unipotent torsion.
If the Weyl group is nontrivial then it has even order;
 in particular, the group $G$ is not of degenerate type in this case.
\end{corollary}

\section*{Acknowledgments}

The authors thank Tuna \Altinel for Lemma \ref{BBC_pi} and
his assistance on various incarnations of the Weyl group argument,
as well as numerous helpful conversations and correction.
We also thank Alexandre Borovik, Adrian Deloro, and Eric Jaligot
for encouragement and enlightening conversations.
The first author gratefully acknowledges the hospitality and support
of Universit\'e Claude Bernard (Lyon-I), University of Manchester,
and Emily Su.

\small
\bibliographystyle{alpha} 
\bibliography{burdges,fMr}

\end{document}